\input amstex
\documentstyle{amsppt}
\magnification=1200
\pageheight{9.5 true in}
\NoBlackBoxes

\font\bg=cmr17

\def\th{\theta}
\def\ft{)_{\theta-1}}
\def\fin{)_\infty}
\def\dq{\,d_q}
\def\Vt{V^{\th}}
\def\ve{\,|\,}
\def\b{\beta}
\def\G{\Gamma}
\def\tht{\thetag}
\def\db{\,d\b(y\ve x)}
\def\const{\operatorname{const}}
\def\Ls{\Lambda^*}
\def\Vs{V^{*\th}}
\def\dbs{\,d\b^*(y\ve x)}
\def\x{x^*}
\def\y{y^*}
\def\xt#1{x_{#1}\,t^{1-#1}}
\def\l{\lambda}
\def\Psm{P^*_\mu}
\def\Psmp{P^*_{\mu'}}
\def\ps{\prec^*}
\def\Is{I^*(\mu,n)}
\def\pmu{{}^\backprime\mu}
\def\f{f_{\mu,\nu}(d)}
\def\fe{f_{\mu,\eta}(d)}
\def\p{\psi_{\mu,\nu}}
\def\i{{}^{(i)}}
\def\ip{{}^{(i+1)}}
\def\pk{p^*_{k}}
\def\os{\omega^*_{q,t}}
\def\Ugln{\Cal U(\frak{gl}(n))}
\def\T{\sssize T}
\def\Sm{{\Bbb S}_\mu}
\def\tr{\operatorname{tr}}
\def\Mat{\operatorname{Mat}}
\def\Q{{\Bbb Q}}

\def\fp#1#2{\langle#1\rangle_{#2}}

\def\mn{\mu/\nu}
\def\fb{\overline{f}}

\leftheadtext{Andrei Okounkov}
\rightheadtext{(Shifted) Macdonald polynomials}

\topmatter
\title\nofrills
\bg
(Shifted) Macdonald Polynomials:\\
q-Integral Representation\\
and Combinatorial Formula
\endtitle
\author
Andrei Okounkov
\endauthor
\affil
Institute for Advanced Study, Princeton
\endaffil
\address
Department of Mathematics, University of Chicago,  5734 South
University Avenue, Chicago, IL 60637-1546
\endaddress
\email
okounkov\@math.uchicago.edu
\endemail
\abstract
We extend some results about {\it shifted Schur
functions\/}  to the general context of
{\it shifted Macdonald polynomials}. We strengthen  
some  theorems of F.~Knop and S.~Sahi and 
give two explicit formulas for these polynomials:
a $q$-integral representation and a combinatorial
formula. Our main tool is a $q$-integral
representation for ordinary Macdonald polynomial.
We also discuss duality for shifted Macdonald
polynomials and Jack degeneration of these polynomials.
\endabstract
\thanks
The author is grateful to the Institute for Advanced
Study for hospitality and to NSF for financial support under
grant DMS 9304580.
\endthanks
\keywords Macdonald polynomials, shifted Macdonald polynomials,
combinatorial formula, $q$-integral representation
\endkeywords
\subjclass 05E05
\endsubjclass
\endtopmatter

\document

\head
\S 1 Introduction
\endhead

The orthogonality of Schur functions
$$
(s_\mu,s_\lambda)=\delta_{\mu\l}\,,
$$
is the orthogonality relation for characters
of the unitary group $U(n)$. The orthogonality
of characters means that a character (as a function
on the group) vanishes in all but one irreducible
representation.

There is a very remarkable basis in the center of the 
universal enveloping algebra $\Ugln$, which is similar
to the basis of characters in the algebra of central
functions on $U(n)$. Many properties of this basis
can be found in \cite{OO,Ok}. The elements of this basis
 are indexed by partitions
$\mu$ with at most $n$ parts; they are denoted by $\Sm$
and called {\it quantum immanants}. The element $\Sm$
has degree $|\mu|$ 
and vanishes in as many
irreducible representations, as possible, namely 
$$
s^*_\mu(\l)=0, \quad\text{unless}\quad \mu\subseteq\l
$$
where $s^*_\mu(\l)$ is the eigenvalue of $\Sm$
 in the representation with highest weight
$\l$. The function $s^*_\mu(\l)$ is called the {\it shifted Schur function}
for it is a polynomial in $\l$ with highest term $s_\mu(\l)$.
This polynomial is symmetric in variables
$$
\l_i - i, \quad i=1,2,\dots \,,
$$
and is also stable, which means it does not change if we
add some zeroes to $\l$. The number $s^*_\mu(\l)$ has a
remarkable combinatorial interpretation: up to some
simple factors, it counts standard tableaux on $\l/\mu$.

There is the following formula for $\Sm$ in terms of
the generators $E_{ij}$ of $\Ugln$. Fix any standard
tableau $T$ on $\mu$. Let $c_{\T}(i)$ denote the content
of the $i$-th square in $T$. Then we have (see \cite{Ok1},
and also \cite{N,Ok2})
$$
\Sm=\tr\left(
(E-c_{\T}(1))\otimes \dots \otimes(E-c_{\T}(|\mu|)) \cdot P_T
\right)\,, \tag 1.1 
$$
where $E$ is the following matrix with entries in $\Ugln$
$$
E=\left(E_{ij}\right)_{ij} \in \Mat(n)\otimes\Ugln\,,
$$
$P_T$ is the orthogonal projection on the Young basis
vector indexed by $T$
$$
P_T \in \Q S(|\mu|)\,,
$$ 
and we use the standard representation
$$
\Q S(|\mu|) \to \Mat(n)^{\otimes|\mu|} \,.
$$
The formula \tht{1.1} corresponds to the following
{\it combinatorial formula} for the shifted Schur
functions (see \cite{OO}, \S11, and also \cite{Ok1}, \S3.7). 
 We call a tableau $T$ on a diagram $\mu$
a {\it reverse tableau} if its entries strictly decrease
down the columns and weakly decrease in the rows. Denote
by  $T(s)$  the entry  of $T$ in the square $s$ and by
$c(s)$ the content of the square $s$. Then we have 
$$
s^*_\mu(x_1,x_2,\dots)=\sum_T \prod_{s\in\mu}
\big(x_{\T(s)}-c(s)\big) \,, \tag 1.2 
$$
where $T$ ranges over all reverse tableau on $\mu$.

The shifted Schur functions have numerous applications to
the finite- and especially infinite-dimensional
representation theory. Their $(q,t)$-analogs, with
which we deal in this paper, inherit most of their power.

In the general $(q,t)$-situation  the center of $\Ugln$
gets replaced by the commutative algebra generated by
Macdonald  $q$-difference operators \cite{M}.
The eigenfunctions of this commutative algebra are
Macdonald polynomials $P_\l(q,t)$, which replace Schur functions.
The eigenvalue of a central
element in the representation with highest weight $\l$
becomes the eigenvalue of a $q$-difference operator on 
the Macdonald polynomial $P_\l(q,t)$. It is known (see, for
example \cite{EK})
that the algebra generated by Macdonald  operators can be
naturally identified with the center of the $q$-deformed $\Ugln$.

The eigenvalue of a Macdonald operator on $P_\l(q,t)$
is known to be a polynomial in $q^{\l_i}$ which
is symmetric in variables 
$$
q^{\l_i} t^{-i}\,.
$$
Therefore the natural $(q,t)$-analog of the shifted
Schur function should be a polynomial
$$
\Psm(x)\,,
$$
of degree
$$
\deg \Psm(x)=|\mu|\,,
$$
which is symmetric in variables
$$
x_i\, t^{-i}, \quad i=1,2,\dots \,,
$$
and satisfies the following vanishing condition
$$
P^*_\mu(q^\l)=0, \quad\text{unless}\quad \mu\subseteq\l \,, \tag 1.3
$$
and $\Psm(q^\mu)\ne 0$. 
It is an overdetermined system of linear conditions on
$\Psm(x)$. It is easy to see that these
polynomials are unique within a scalar factor provided
they exist. We shall specify the normalization in section 4.
We call these polynomials {\it shifted
Macdonald polynomials}. 

It is well known that the number $q^{j-1}t^{1-i}$ is the $(q,t)$-analog
of the content of the a square $s=(i,j)$. Together with G. Olshanski
we have conjectured (unpublished) the following analog of the
formula \tht{1.2}. By
$$
a'(s)=j-1\,,\quad l'(s)=i-1
$$
denote the {\it arm-colength\/} and the
{\it leg-colength\/} of the square $s=(i,j)$. Then 
$$
\Psm(x;q,t)=\sum_T \psi_{\T}(q,t) \prod_{s\in\mu} t^{1-\T(s)}
\left(x_{\T(s)}-q^{a'(s)} t^{-l'(s)}
\right) \,, \tag 1.4
$$
where the sum is over all reverse tableau on $\mu$ with entries
in $\{1,2,\dots\}$  and $\psi_{\T}(q,t)$ is the same $(q,t)$-weight of a tableau
which enters the combinatorial formula for ordinary Macdonald polynomials
(see \cite{M},\S VI.7)
$$
P_\mu(x;q,t) = \sum_T \psi_{\T}(q,t) \prod_{s\in\mu} x_{\T(s)} \,.  \tag 1.5
$$
The coefficients $\psi_T(q,t)$ are rational functions of $q$ and $t$.

However, only in the present paper we obtain a proof of this formula. 

First theorems  about the polynomials $\Psm(x)$ were obtained
by F.~Knop and S.~Sahi \cite{KS,K,S}. In particular, they identified
the highest degree term of the inhomogeneous polynomial $\Psm(x)$. 
Other fundamental properties (such as {\it integrality}) were 
also established. Their approach was based on 
$q$-difference equations for polynomials $\Psm(x)$. 

The combinatorial formula \tht{1.4} implies both the highest
degree term identification \cite{K,S} and the {\it extra vanishing}
property \tht{1.3} proved in \cite{K}. It is easy to see that all conditions
from the definition of $\Psm(x)$ are obvious in \tht{1.4} except
for the  symmetry in $x_i t^{-i}$. We were not able to find any direct proof 
of this symmetry and shall give a quite indirect proof. 

We shall use an approach based
on a {\it $q$-integral representation} for the polynomial $\Psm(x)$.
It is independent of the difference equations approach
and generalizes the {\it coherence} property of 
quantum immanants (see \cite{OO}, \S10 and also \cite{Ok1}, \S5.1).   

Our main technical tool is a $q$-integral representation for
Macdonald polynomials (see Theorem I below). Since
Macdonald polynomials satisfy a 
$q$-difference equation it is natural to expect
that 
$$
P_\mu(x_1,\dots,x_n)  
$$
can be written as a multiple $q$-integral where $x_i$ occur
as limits of integration. Indeed, $P_\mu(x_1,\dots,x_n)$  can be written as
such a multiple integral of the Macdonald polynomial
$$
P_\mu(y_1,\dots,y_{n-1})
$$
in a smaller set of variables with respect to a 
multivariate analog of the symmetric beta measure.
In the simplest case $n=2$ this integral reduces to
a particular case of the $q$-analog of the beta integral
studied in \cite{AA,AV}. In the Schur function
case this integral can be evaluated explicitly
and gives the determinant ratio formula. The integral
representation allows to characterize $P_\mu(x)$ (and also $\Psm(x)$)
as eigenfunctions of commuting integral operators, see
remark 4.7.

There are important $q$-integral formulas involving 
Macdonald polynomials due to Kadell (see \cite{M}, example
VI.9.3 and also \cite{Ka})
as well as integral representation of Macdonald
polynomials via ordinary contour integrals
(see \cite{AOS} and references therein).

The crucial property
of our integral is that the domain of integration  is of the form
$$
\int_{x_2}^{x_1} \dq y_1 \dots \int_{x_n}^{x_{n-1}} \dq y_{n-1} 
\, \big(\,\dots\,\big)\,. 
$$
This allows to obtain a $q$-integral representation
for the polynomial $\Psm$ just by a minor modification
of the integral, see section 4. In other words, the
orthogonal polynomials $P_\mu$ and the Newton interpolation
polynomials $\Psm$ have essentially the same integral
representation.

In a sense the relationship between $P_\mu$ and $\Psm$
is even closer than between $s_\mu$ and $s^*_\mu$. 
One explanation for this is that the finite difference
calculus (which is suitable to $s^*_\mu$) unifies with
the ordinary calculus (which is suitable to $s_\mu$)
in the $q$-difference calculus.

The lower degree terms of the inhomogeneous
polynomial $\Psm$ carry some important additional
information and make some properties of $\Psm$ look
even more natural than the corresponding properties
of the ordinary Macdonald polynomials. For example,  
the {\it duality} for shifted Macdonald polynomials has
a clear combinatorial interpretation  (see section 6) which gives a
new interpretation of the duality for Macdonald
polynomials.

In section 7 we consider shifted Jack polynomials.
This degeneration was also considered by F.~Knop and
S.~Sahi in \cite{KS}. 
Even if one is interested in Jack polynomials
only it proves to be easier to work with general Macdonald
polynomials and then let $q\to 1$ in very final formulas.

It would be 
 very interesting to find a $(q,t)$-analog of the formula \tht{1.1}.

I am grateful to I.~Cherednik, S.~Sahi and especially to G.~Olshanski
for many helpful discussions. I would like to thank R.~Askey, 
A.~N.~Kirillov, K.~Mimachi, M.~Noumi, M.~Wakayma for their
remarks on the preliminary version of this paper (q-alg 9605013).
In particular, K.~Mimachi gave me a copy of the 
preprint \cite{MN}. 

Since the completion of this paper  binomial type
formulas for $\Psm$ and also for ordinary Macdonald 
polynomials were obtained (see \cite{Ok3} and also
\cite{OO3}). In particular, they provide
new and, perhaps, more natural ways to prove the
$q$-integral representations obtained here.    
Among other application let us mention the papers
\cite{KOO} and \cite{Ok5}.

The analogs of $\Psm$ for the root system of type $BC_n$
were considered in \cite{Ok4}. Those polynomials have
essentially all properties of $\Psm$ except for the 
difference equations. 

\head
\S2 Notations
\endhead

All $q$-shifted factorials 
in this paper will be with the same base $q$, for example
$$
(a\fin=(1-a)(1-q a)(1-q^2 a)\dots\,.
$$
This product converges if $|q|<1$.
Put
$$
(a)_\th=\frac{(a\fin}{(q^\th a\fin} \,.
$$
The numbers $q$ and $t=q^\th$ will be the two parameters of the
Macdonald polynomials $P_\mu$. With this notation
$$
(a)_\th=\frac{(a\fin}{(t a\fin} \,.
$$
We shall need also another $q$-shifted power 
$$
\fp{a}{r}=(a-1)(a-q)\dots(a-q^{r-1}),\quad r=0,1,2,\dots\,,
$$
which is a $q$-analog of the falling factorial power of $a$.

Recall the definition of the
$q$-integral
$$
\int_b^a f(y)\dq y = \int_0^a f(y)\dq y - \int_0^b f(y)\dq y\,,
$$
where
$$
\int_0^a f(y)\dq y=a(1-q)\sum_{i=0}^\infty
f(a q^i) q^i \,.
$$
We have
$$
\int_0^a y^r\dq y = \frac1{[r+1]} a^{r+1} \,, \tag 2.1
$$
where
$$
[r]=\frac{1-q^r}{1-q}
$$
is the $q$-analog of the number $r$. Recall also the following
$q$-analog of the beta function integral \cite{GR,1.11.7}
$$
\int_0^1 y^{a-1} (q y)_{b-1} \dq y = B_q(a,b)\,, \tag 2.2
$$
where $\Re a>0$, $b\ne0,-1,-2,\dots$,
$$
B_q(a,b)=\frac{\G_q(a) \G_q(b)}{\G_q(a+b)}\,,
$$
and, finally,
$$
\G_q(a)=(1-q)^{1-a} (q)_{a-1} \,.
$$

Given two vectors
$$
x=(x_1,\dots,x_n),\quad y=(y_1,\dots,y_{n-1})\,,
$$
where $n=2,3,\dots$, write
$$
y\prec x
$$
if 
$$
y_i\in [x_i,x_{i+1}],\quad i=1,\dots,n-1\,.
$$
Denote by
$$
V(x)=\prod_{i<j} (x_i-x_j)
$$
the Vandermonde determinant in variables $x_1,\dots,x_n$. Put
$$
\Vt(x)=V(x) \prod_{i\ne j} (q x_i/x_j\ft \,.
$$
We will integrate symmetric polynomials in $y_1,\dots,y_{n-1}$
over the domain $y\prec x$ with respect to the following
measure
$$
d\b(y\ve x)= V(y) \prod_{i,j} (q y_i/x_j\ft \dq y\,, \tag 2.3
$$
which is a multivariate analog of the beta measure in \tht{2.2}.
Here
$$
d_q y=d_q y_1 \dots d_q y_{n-1}\,.
$$
With the Macdonald notation
$$
\Pi(x,y;q,t)=\prod\frac{(t x_i y_j\fin}{(x_i y_j\fin}
$$
we have
$$
d\b(y\ve x)= V(y)\, \Pi(1/x,t y;q,q/t) \dq y\,, \tag 2.4
$$

By $\Lambda$ denote the algebra of symmetric functions
with coefficients in rational functions in $q$ and $t$.
By $\Lambda(n)$ denote the algebra on symmetric functions
in $n$ variables with the same coefficients.

\head
\S3\,$q$-Integral representation of Macdonald polynomials.
\endhead

Given a partition $\mu$ with the number of part $\ell(\mu)\le n$
put
$$
C(\mu,n)=\prod B_q(\mu_i+(n-i)\th,\th)\,. \tag 3.1
$$
In this section we will prove the following theorem

\proclaim{Theorem I}
Suppose $\ell(\mu)< n$, then
$$
\frac1{\Vt(x)}\int_{y\prec x} P_\mu(y)\db  = C(\mu,n) P_\mu(x) \,.
\tag 3.2
$$
\endproclaim

On the one hand this formula expresses the Macdonald 
polynomial $P_\mu(x)$ in terms of the Macdonald polynomial $P_\mu(y)$
in the smaller set of variables. By stability of
Macdonald polynomials it suffices to know $P_\mu$
in $|\mu|$ variables. Therefore \tht{3.2} gives a
$q$-integral representation of $P_\mu$. This formula is
a generalization of the determinant ratio formula for
Schur functions (see below). 

On the other hand, since $P_\mu$ form a basis in the space
of symmetric polynomials, \tht{3.2} tells how to integrate
symmetric polynomials over the domain $y\prec x$ with
respect to the measure \tht{2.3}.

In the limit $q\to1$ the formula \tht{3.2} becomes the
integral representation for Jack polynomials found by
G.~Olshanski \cite{Ol}.

Both sides of \tht{3.2} are analytic functions of $q$ and $t$ in
the polydisc
$$
|q|<1,\,|t|<1\,.
$$
More precisely, if
$$
|q|,|t|<\delta<1\,,
$$
then we have to assume, for example, that
$$
\delta^{1/2}<|x_i|<\delta^{-1/2},\quad i=1,\dots,n
$$
and that
$$
x_i\ne x_j,\quad i\ne j
$$
in order to avoid zero factors.

Therefore it suffices to prove the equality \tht{3.2}
for 
$$
t=q^\th, \quad \th=2,3,\dots\,. \tag 3.3
$$
Indeed, an analytic function, which vanishes on \tht{3.3}
should vanish on all lines 
$$
q=\const\ne 0
$$
by the 1-dimensional uniqueness theorem. 

\proclaim{Lemma 3.1} Suppose $f(y_1)$ is a polynomial. Then
$$
\frac 1{x_1-x_2}\int_{x_2}^{x_1} f(y_1)\dq y_1 \tag 3.4
$$
is a symmetric polynomial in $x_1$ and $x_2$.
\endproclaim
\demo{Proof} Follows from \tht{2.1}. \qed
\enddemo

The following computation is from \cite{OO}, \S10.

\proclaim{Lemma 3.2} Suppose $f(y)$ is a symmetric polynomial. Then
$$
\int_{y\prec x} f(y) V(y) \dq y \tag 3.5
$$
is a skew-symmetric polynomial in $x$.
\endproclaim
\demo{Proof} Since $f(y) V(y)$ is skew-symmetric, it 
is a linear combination of the following determinants
$$
\det
\left(
\matrix
y_1^{\nu_1}&\hdots& y_1^{\nu_{n-1}}\\
\vdots&&\vdots\\
y_{n-1}^{\nu_1}&\hdots& y_{n-1}^{\nu_{n-1}}
\endmatrix
\right)
$$
for some numbers $\nu_1>\nu_2>\dots$. Integrating we get
up to a constant factor
$$
\det
\left(
\matrix
x_1^{\nu_1+1}-x_2^{\nu_1+1}&\hdots&
x_1^{\nu_{n-1}+1}-x_2^{\nu_{n-1}+1}\\
\vdots&&\vdots\\
x_{n-1}^{\nu_1+1}-x_n^{\nu_1+1}&\hdots&
x_{n-1}^{\nu_{n-1}+1}-x_n^{\nu_{n-1}+1}
\endmatrix
\right)\,.
$$
Which equals
$$
\det
\left(
\matrix
x_1^{\nu_1+1}&\hdots& x_1^{\nu_{n-1}+1}&1\\
x_2^{\nu_1+1}&\hdots& x_2^{\nu_{n-1}+1}&1\\
\vdots&&\vdots&\vdots\\
x_n^{\nu_1+1}&\hdots&x_n^{\nu_{n-1}+1}&1
\endmatrix
\right)\,. \tag 3.6
$$
To see this subtract in \tht{3.6} the second line from
the first one, then the third line from the second one
and so on. Clearly, \tht{3.6} is skew-symmetric in $x$.
\qed
\enddemo

\remark {\bf Remark 3.3} In the Schur functions case $q=t$
the computation from the previous lemma gives
in fact an explicit computation of the integral 
in the LHS of \tht{3.2} which shows that \tht{3.2} is
a generalization of the determinant ratio formula
for Schur functions.
\endremark

Denote by $I(\mu,n)$ the LHS of \tht{3.2}.

\proclaim{Proposition 3.4} Suppose $f(y)$ is 
symmetric polynomial. Then
$$
\frac1{\Vt(x)}\int_{y\prec x} f(y)\db \tag 3.7
$$
is a symmetric polynomials in $x$ of degree $\deg f$.
\endproclaim
\demo{Proof} Multiply both the integrand and denominator
in \tht{3.7} by 
$$
\left(\prod x_i\right)^{(n-1)(\th-1)} 
$$
so that they become polynomials. The integrand becomes
a skew-symmetric polynomial in $y$ with coefficients
in symmetric polynomials in $x$, therefore by lemma 3.2
the integral is a skew-symmetric polynomial in $x$.
Denote this polynomial by $J$.

By lemma 3.1 the polynomial $J$ is divisible by $(x_1-x_2)$.

Since the integrand vanishes at the points 
$$
y_1=x_1/q,\dots,x_1/q^{\theta-1}\,,
$$
as well as at the points 
$$
y_1=x_2/q,\dots,x_2/q^{\theta-1}\,,
$$
we can replace the integration
$$
\int_{x_2}^{x_1} \dq y_1
$$
by integration
$$
\int_{x_2}^{x_1/q^s} \dq y_1,\quad s=1,\dots,\theta-1,
$$
or by
$$
\int_{x_2/q^s}^{x_1} \dq y_1,\quad s=1,\dots,\theta-1\,.
$$
Therefore $J$ is divisible also by
$$
x_1-q^s x_2, \quad s=1,\dots,\theta-1,
$$
and 
$$
q^s x_1- x_2, \quad s=1,\dots,\theta-1.
$$

Since $J$ is skew-symmetric in $x$ it is divisible
by 
$$
V(x) \prod_{i\ne j} \prod_{s=1}^{\theta-1} (x_i-q^s x_j) \,.
$$ 
Therefore \tht{3.7} is a symmetric polynomial
in $x$. It is clear that its degree equals 
$$
\deg f+(n-1)(n-2)/2+(n-1)-n(n-1)/2=\deg f\,,
$$
where the summands come from $f(y)$, $V(y)$,
integration, and $\Vt(x)$ respectively. This
concludes the proof. \qed
\enddemo

Denote by $I(\mu,n)$ the LHS of \tht{3.2}.
We want to apply the Macdonald operator
$$
D=\sum_{i=1}^n \prod_{j\ne i} \frac{t x_i-x_j}{x_i-x_j} \,T_{q,x_i} \tag 3.8
$$
to $I(\mu,n)$. Here
$$
\left[T_{q,x_i} f\right](x_1,\dots,x_i,\dots,x_n)=
f(x_1,\dots,q x_i,\dots,x_n) \,.
$$
The operator \tht{3.8} will be also denoted by $D({q,t})$ and by $D_x({q,t})$
when it should be stressed that it acts on variables $x$.
The Macdonald polynomials are eigenfunctions of this operator
\cite{M}
$$
D\, P_\mu=\left(\sum q^{\mu_i} t^{n-i}\right) P_\mu\,. \tag 3.9
$$

We shall need the three following general lemmas about the
operator $D$. For these lemmas we do not need any special
assumptions about $q$ and $t$.

\proclaim{Lemma 3.5} For all $q$ and $t$ we have the following
commutation relation
$$
D({1/q,1/t}) \frac1{\Vt(x)}=\frac{(q/t)^{(n-1)}}{\Vt(x)} D({1/q,t/q})\,.
\tag 3.10
$$
\endproclaim
\demo{Proof} Direct computation. \qed
\enddemo

\proclaim{Lemma 3.6} Put
$$
D_{1/x}(q,t)=\sum_i \prod_{j\ne i} \frac{t/x_i-1/x_j}
{1/x_i-1/x_j} \,T_{q,1/x_i}\,. 
$$
then
$$
D_{1/x}(q,t)=t^{n-1} D_x(1/q,1/t) \,. \tag 3.11
$$
\endproclaim
\demo{Proof} Direct computation. \qed
\enddemo

By definition \tht{3.8} the operator $D$ acts on symmetric
polynomials in $n$ variables. This action indeed depends
on $n$; in other words it is not compatible with
the restriction homomorphisms
$$
\Lambda(n) \to \Lambda(n-1)\,. \tag 3.12
$$

In the next lemma we shall deal with two finite
sets of variables $x_1,\dots,x_n$ and $y_1,\dots,y_m$
and we suppose for simplicity that $n\ge m$. Denote
by $D_x$ and $D_y$ the operators $D$ in variables
$x$ and $y$ respectively. Put
$$
\Pi_{n,m}=\Pi(x,y;q,t)\,.
$$

\proclaim{Lemma 3.7} 
$$
D_x \Pi_{n,m} = \left( t^{n-m} D_y + [n-m]_t\right) \Pi_{n,m}\,, 
\tag 3.14
$$
where $[n-m]_t=(1-t^{n-m})/(1-t)$.
\endproclaim
\demo{Proof} 
As explained in \cite{M, VI.4} (this is clear from \tht{3.9})
the following modification of the operator \tht{3.8}
$$
E=t^{-n} D - \sum_1^n t^{-i} \tag 3.15
$$
is compatible with homomorphisms \tht{3.12} and therefore
defines an operator
$$
\Lambda\to\Lambda\,,
$$
which is self-adjoint. By \cite{M, VI.2.13}
$$
E_x \Pi(x,y;q,t) = E_y \Pi(x,y;q,t) \,.
$$
Therefore
$$
\align
D_x \Pi_{n,m} &= 
t^n \left(E_x + \sum_1^n t^{-i}\right) \Pi_{n,m} \\
&= t^n \left(E_y + \sum_1^n t^{-i}\right) \Pi_{n,m} \\ 
&= t^n \left(t^{-m} D_y - \sum_1^m t^{-i}+ \sum_1^n t^{-i}\right) \Pi_{n,m} \\
&= \left(t^{n-m} D_y +[n-m]_t\right) \Pi_{n,m} \,. \qed
\endalign
$$
\enddemo

In particular,
$$
D_x \Pi_{n,n-1} = \left(t D_y + 1\right) \Pi_{n,n-1}\,.\tag 3.16
$$

\proclaim{Proposition 3.8}
$$
D(1/q,1/t) I(\mu,n) = \left(\sum_i q^{-\mu_i} t^{i-n} \right) I(\mu,n)\,.
\tag 3.17
$$
\endproclaim
\demo{Proof}
Recall that we consider the case $\theta=2,3,\dots$.
Put
$$
I_0 = \int_{y\prec x} P_\mu(y)\db \,.
$$
By \tht{3.10} we have to prove that
$$
(q/t)^{n-1} D(1/q,t/q) I_0 = \left(\sum_i q^{-\mu_i} t^{i-n} \right) I_0 \,.
$$
By \tht{3.11} it is equivalent to
$$
D_{1/x}({q,q/t}) I_0 = \left(\sum_i q^{-\mu_i} t^{i-n} \right) I_0 \,.
$$

We have
$$
T_{q,1/x_i} I_0 = 
\int_{x_2}^{x_1} \!\dots \!\int_{x_i/q}^{x_{i-1}} 
\!\int_{x_{i+1}}^{x_i/q}\!\dots\!\int_{x_n}^{x_{n-1}}
T_{q,1/x_i}  P_\mu(y)\db \,. \tag 3.18
$$
Since the integrand in \tht{3.18} vanishes if
$$
y_j,y_{j-1}=x_{j}/q, \quad j\ne i,
$$
we can rewrite \tht{3.18} as follows
$$
T_{q,1/x_i} I_0 =\int_{y\prec x/q} T_{q,1/x_i} P_\mu(y)\db \,. 
$$
Therefore
$$
D_{1/x}(q,q/t) I_0 =\int_{y\prec x/q} D_{1/x}(q,q/t) P_\mu(y)\db \,. 
\tag 3.19
$$
Now by \tht{3.16} 
$$
D_{1/x}(q,q/t) \Pi(1/x,t y;q,q/t) = 
\left(\frac qt D_y(q,q/t) + 1\right) \Pi(1/x,t y;q;q/t) \,. \tag 3.20
$$
By \tht{3.20} and \tht{2.4}, the integral \tht{3.18} equals
$$
\int_{y\prec x/q} P_\mu(y) V(y) 
\left[
\left(\frac qt D_y(q,q/t) + 1\right) \Pi(1/x,t y;q;q/t)
\right] \dq y \,.\tag 3.21
$$ 
Consider one summand in \tht{3.21}
$$
\frac qt \int_{y\prec x/q} P_\mu(y) V(y) 
\left[ \prod_{j\ne i} 
\frac{q y_i/t-y_j}{y_i-y_j}\, T_{q,y_i} \,
\Pi(1/x,t y;q;q/t) \right] \dq y \,,\tag 3.22
$$ 
where $i=1,\dots,n-1$. Replace $q y_i$ by a new variable.
To simplify notation denote this new variable by $y_i$. Then
\tht{3.22} becomes
$$
\frac 1t \int_{x_2/q}^{x_1/q} \!\dots \!\int_{x_{i+1}}^{x_{i}} 
\!\dots\!\int_{x_n/q}^{x_{n-1}/q} 
\left[ \prod_{j\ne i} 
\frac{y_i/t-y_j}{y_i-y_j} \,T_{1/q,y_i} \,
P_\mu(y)\right] \db \,. \tag 3.23
$$
Since the beta measure vanishes if
$$
y_j=x_j/q,x_{j+1}/q, \tag 3.24
$$
the integral \tht{3.23} equals
$$
\frac 1t \int_{y\prec x} 
\left[ \prod_{j\ne i} 
\frac{y_i/t-y_j}{y_i-y_j} \,T_{1/q,y_i}\,
P_\mu(y)\right] \db \,. \tag 3.25
$$
By the same vanishing of the beta measure for \tht{3.24}
$$
\int_{y\prec x/q} P_\mu(y) V(y) 
\Pi(1/x,t y;q;q/t)
\dq y = 
\int_{y\prec x} P_\mu(y) V(y) 
\Pi(1/x,t y;q;q/t)
\dq y \,.\tag 3.26
$$ 
Therefore the integral \tht{3.21}
equals
$$
I_0+\frac 1t \int_{y\prec x} 
\left[ D_y(1/q,1/t) P_\mu(y)\right] \db \,, \tag 3.27
$$
It is well known (and follows, for example, from the
formula for the scalar product) that
$$
P_\mu(x;{q,t})=P_\mu(x;1/q,1/t)\,. \tag 3.28
$$
By \tht{3.9} and \tht{3.28} 
$$
D_y(1/q,1/t) P_\mu(y) = \left(\sum_1^{n-1} q^{-\mu_i} t^{i-n+1}
\right) P_\mu(y)\,.
$$
Hence \tht{3.27} equals
$$
\left(\sum_1^{n} q^{-\mu_i} t^{i-n}
\right) I_0\,,
$$
as desired. \qed
\enddemo

\proclaim{Corollary 3.9}
$I(\mu,n)$ equals the Macdonald polynomial $P_\mu(x)$ up to
a scalar factor.
\endproclaim

We calculate this factor in the
next proposition. It will conclude the proof of the
theorem.

Consider the highest monomial in  of $I(\mu,n)$
with respect to the lexicographic ordering of monomials

\proclaim{Proposition 3.10}
The highest monomial in $I(\mu,n)$ equals
$$
C(\mu,n)\, \prod  x_i^{\mu_i} \,. \tag 3.29
$$
\endproclaim
\demo{Proof} Multiply both the integrand and denominator
in $I(\mu,n)$ by 
$$
\prod (-1)^{(\th-1)(n-i)} q^{-(n-i)\th(\th-1)/2} x_i^{(n-1)\th}\,. 
$$
Then the highest monomial in the denominator equals
$$
\prod_i x_i^{(2\th-1)(n-i)}\,. 
$$

Calculate the highest term of the integrand. We have to
give $x_i$ and $y_i$ the same priority because we obtain $x_i$
integrating $y_i$. Therefore the highest term of the integrand
is
$$
\align
&\prod_i x_i^{(n-i-1)(\th-1)}
y_i^{\mu_i+(n-i-1)+(n-i)(\th-1)}
 (x_i-q y_i)\dots(x_i-q^{\th-1}y_i) =\\
&\prod_i x_i^{(n-i-1)(\th-1)}
y_i^{\mu_i+(n-i)\th-1}
(x_i-q y_i)\dots(x_i-q^{\th-1}y_i) \,.
\endalign
$$

Now calculate the highest term of the integral. The terms
which come from the lower limits are negligible, therefore
the highest term of the integral equals
$$
\align
&\prod_i x_i^{(n-i-1)(\th-1)} \int_0^{x_i} 
y_i^{\mu_i+(n-i)\th-1}
(x_i-q y_i)\dots(x_i-q^{\th-1}y_i) \dq y_i \\
&=\prod_i x_i^{(n-i-1)(\th-1)+\mu_i+(n-i)\th+\th-1} \int_0^{1} 
z^{\mu_i+(n-i)\th-1}
(1-q z)\dots(1-q^{\th-1} z) \dq z \\
&=\prod_i 
B_q (\mu_i+(n-i)\th,\th)\,x_i^{\mu_i+(n-i)(2\th-1)}\\
&=C(\mu,n) \, \prod_i x_i^{\mu_i+(n-i)(2\th-1)}\,, 
\endalign
$$
where we used change of variables $y_i=z x_i$, beta
function integral \tht{2.2} and definition \tht{3.1}. 
Therefore the highest term of the ratio equals
$$
C(\mu,n)\,\prod_i x_i^{\mu_i}\,.\qed
$$ 
\enddemo

\example{Example 3.11}
The integral representation \tht{3.2} gives one
more way to calculate the special value (see \cite{M}, VI.6.17)
$$
P_\mu(1,t,t^2,\dots,t^{n-1})\,. \tag 3.30
$$
What is special about this value is that for
$$
\theta=1,2,\dots
$$
only one summand, corresponding to the point
$$
y_i=t^{i-1}
$$
does not vanish in \tht{3.2}.
\endexample

\head
\S4\,$q$-Integral representation for shifted Macdonald 
polynomials
\endhead

Denote by $\Ls(n)$ the algebra of polynomials
$f(x_1,\dots,x_n)$ which are symmetric in new
variables 
$$
x'_i=\xt{i} \,.
$$
We shall call such polynomials {\it shifted
symmetric}. The algebra $\Ls(n)$ is filtered
by degree of polynomials.
Denote by $\Ls$ the projective limit of 
filtered algebras
$\Ls(n)$ with respect to the homomorphisms
$$
\align
\Ls(n)&\to\Ls(n-1)\,,\\
f(x_1,\dots,x_n)&\mapsto f(x_1,\dots,x_{n-1},1)\,.
\endalign
$$
The algebra $\Ls$ can be naturally identified
with the algebra of those polynomials in 
Macdonald commuting difference operators which
are stable (that is compatible with homomorphisms \tht{3.12}).
An example of such operator is the operator $E$
from \tht{3.15}. 

The graded algebra corresponding to the filtered 
algebra $\Ls$ can be naturally identified with
the algebra of symmetric functions in 
variables $x'_i$. The algebra $\Ls$
is generated, for example, by the following analogs
of the power sums
$$
p^*_k(x)=\sum_i (x_i^k-1)t^{k(1-i)},\quad k=1,2,\dots \,. \tag 4.1
$$

Given a partition $\mu$ put
$$
n(\mu)=\sum_i (i-1) \mu_i = \sum_j \mu'_j(\mu'_j-1)/2 \,.
$$
Recall that for each square $s=(i,j)\in\mu$ the numbers
$$
\alignat2
&a(s)=\mu_i-j,&\qquad &a'(s)=j-1,\\
&l(s)=\mu'_j-i,&\qquad &l'(s)=i-1,
\endalignat
$$
are called arm-length, arm-colength, leg-length, and
leg-colength respectively. By definition, put
$$
H(\mu)=t^{-2n(\mu)} q^{n(\mu')} \, \prod_{s\in\mu} 
\left(q^{a(s)+1} t^{l(s)} -1\right)\,.
$$
This number will play the same role as the hook-length
product played in \cite{OO}.

Suppose $\ell(\mu)\le n$.
By $\Psm$ denote the element of $\Ls(n)$
which satisfies the two following conditions
$$
\align
\deg \Psm &= |\mu|\,, \tag 4.2\\
\Psm(q^\l)&=H(\mu)\,\delta_{\l,\mu}\,, \quad |\l|\le|\mu|,
\ell(\l)\le n\,, \tag 4.3
\endalign
$$
where 
$q^\l=(q^{\l_1},\dots,q^{\l_n})$. 
Here $H(\mu)$ is just normalization constant and is
introduced for convenience. Note that the condition \tht{4.3}
is weaker than the condition \tht{1.3} in the introduction.
Here we have a square system of linear equations on $\Psm$.
We shall prove \tht{1.3} below in \tht{4.11}.

It is clear that if $\Psm$ exists then it is unique.
The existence of $\Psm$ will follow from an explicit
formula for it. The existence of this polynomial 
was proved by different methods in \cite{K,S}.

In the same way as in \cite{OO}, Ex.~3.5, we have
\proclaim{Proposition 4.1}
The sequence 
$$
\left\{\Psm\in\Ls(n)\right\}_{n\ge\ell(\mu)}
$$
defines an element of $\Ls$.
\endproclaim
\demo{Proof}
The polynomial
$$
\Psm(x_1,\dots,x_{n-1},1)
$$
satisfies all conditions for the shifted Macdonald polynomial
in $(n-1)$ variables provided $n>\ell(\mu)$. \qed
\enddemo 

\proclaim{Proposition 4.2} Suppose $\mu_n>0$ and put
$$
\mu^{-}=(\mu_1-1,\dots,\mu_n-1)\,.
$$
Then
$$
\Psm(x_1,\dots,x_n)=q^{|\mu^{-}|} 
\prod_i \left(\xt{i}-t^{1-n}\right)
\, P^*_{\mu^{-}}(x_1/q,\dots,x_n/q) \,. \tag 4.4
$$
\endproclaim
\demo{Proof} We have to verify that the RHS of \tht{4.4}
satisfies all conditions for $\Psm$. It is clear, that
its degree equals $|\mu|$. Evaluate it at $q^{\mu}$.
We have
$$
q^{|\mu^{-}|} t^{-n(n-1)}
\prod_i \left(q^{\l_i} t^{n-i}-1\right)
\, H(\mu^{-}) = H(\mu) \,. 
$$
Finally, it is easy to see that it vanishes at all points
$q^\l$, where $|\l|\le|\mu|$ and $\ell(\l)\le n$. \qed
\enddemo

In this section we shall obtain a $q$-integral formula
for polynomials $\Psm$, which is a minor modification 
of the formula \tht{3.2}.

We shall consider the following integration 
$$
\int_{y\ps x} \dq y = \int_{x_2}^{q x_1}\dq y_1
\!\dots\!  \int_{x_n}^{q x_{n-1}}\dq y_{n-1} \,.
$$
Introduce new variables
$$
\x_i=x_i\,t^{n-i},\quad \y_i=y_i\,t^{n-1-i}\,,
$$
and put
$$
\align
V^*(x)&=V(\x)\,,\\
\Vs(x)&=\Vt(\x)\,,\\
d\beta^*(y\ve x)&=d\beta(\y\ve\x)\,.
\endalign
$$
We have
\proclaim{Theorem II} Suppose $\ell(\mu)< n$ then
$$
\frac1{\Vs(x)}\int_{y\ps x} \Psm(y)\dbs = t^{|\mu|} C(\mu,n) 
\Psm(x) \,. \tag 4.5
$$
\endproclaim

Together with \tht{4.4} the theorem gives a formula
for polynomials $\Psm$. In the proof we shall assume
that
$$
\th=1,2,\dots\,. \tag 4.6
$$
Indeed, the LHS of \tht{4.5} is analytic in the polydisc
$|q|,|t|<1$ and being a polynomial vanishing at certain
points provided \tht{4.6} must be such a polynomial for
the other values of $q$ and $t$ as well.

We need the two following elementary lemmas.

\proclaim{Lemma 4.3} For any $A$
$$
\frac{(qA/t\ft}{(q/A\ft}=A^{\th-1} q^{-\th(\th-1)/2} \,. \tag 4.7
$$
\endproclaim

\proclaim{Lemma 4.4} For any partition $\mu$ and $\th=1,2,\dots$
$$
\Vs(q^\mu)\ne 0 \,. \tag 4.8
$$
\endproclaim

\demo{Proof of the theorem}
Denote by $\Is$ the integral in the LHS of \tht{4.5}. First show
that it is an element of $\Ls(n)$. In variables $\x,\y$
it can be rewritten as
$$
\frac1{\Vt(\x)}
\int_{\x_2}^{q^{1-\th} \x_1}
\!\dots\!  \int_{\x_n}^{q^{1-\th}\x_{n-1}}
\Psm(y)
\,d\beta(\y\ve\x)\,. \tag 4.9
$$
The polynomial
$$
\Psm(y)
$$
is symmetric in $\y_i$. Because of the vanishing of the
beta measure we can replace
integration
$$
\int_{\x_{i+1}}^{q^{1-\th} \x_i}
$$
by integration
$$
\int_{\x_{i+1}}^{\x_i} \,.
$$
Thus, by proposition 3.4 it is a symmetric polynomial
in $\x_i$, that is an element of $\Ls(n)$.

The main point is the verification of the vanishing
condition. Let us evaluate $\Is$ at
$$
x=q^\l,
$$
where $\l$ is a partition. Then the summands which enter
the $q$-integral correspond to points
$$
y_i=q^{\gamma_i},\quad \l_i\ge\gamma_i\ge\l_{i+1},\quad i=1,\dots,n-1\,.
$$
(Recall that $q<1$ and therefore  $q^{\l_{i+1}}$ is in fact
the upper limit of integration for $y_i$ and $q^{\l_i+1}$ is
the lower one. This explains also the minus sign below in
\tht{4.10}.) Note that $\gamma$ is a partition
$$
\gamma_1\ge\gamma_2\ge\dots\ge\gamma_{n-1}
$$
and $|\gamma|\le|\l|$.

Now suppose that   
$|\l|\le|\mu|$ and $\l\ne\mu$. Then always
$$
\gamma\ne\mu
$$
By definition, $\Psm(y)$
vanishes at all such points! Since by \tht{4.8}
the denominator does not vanish, $\Is$ vanishes.

The only calculation left is evaluation of the integral at  
$$
x=q^\mu\,,
$$
which means
$$
\x_i=q^\mu_i t^{n-i}\,.
$$
In particular $x_n=\x_n=1$.
This evaluation is elementary but quite messy.
 To the end of the proof the
indices $i$ and $j$ will range from $1$ to $n-1$.
Only one summand in the integral does not vanish. This summand
equals the product of the following factors:
$$
H(\mu)
$$
from the value of $\Psm(y)$ at this point,
$$
t^{-(n-1)(n-2)/2} {\prod_{i<j}(\x_i-\x_j)}
$$
from the Vandermonde determinant,
$$
((q/t\ft)^{n-1} {\prod (q\x_i/t\ft} {\prod_{i\ne j} ((q/t) \x_i/\x_j\ft} 
$$
from the other factors in the beta density, and finally
$$
(-1)^{n-1} (1-q)^{n-1} t^{-(n-1)} {\prod \x_i} \tag 4.10
$$
from the q-Lebesgue measure $\dq \y$.

The denominator equals
$$
{\prod_{i<j}(\x_i-\x_j)} {\prod_{i\ne j} (q \x_i/\x_j\ft}
\prod_i \left( (\x_i-1) (q\x_i\ft (q/x_i\ft \right) \,.
$$
By \tht{4.7}
$$
\prod_{i\ne j} \frac{((q/t) \x_i/\x_j\ft}{(q \x_i/\x_j\ft} =
q^{-\th(\th-1)(n-1)(n-2)/2} \,.
$$
Again, by \tht{4.7}
$$
\frac
{(-1)\x_i (1-q) (q/t\ft (q\x_i/t\ft}
{(\x_i-1)(q\x_i\ft (q/\x_i\ft} =
\frac
{(\x_i)^\th q^{\th(\th-1)} (1-q) (q\ft}{(\x_i)_\th}\,,
$$
which equals
$$
(\x_i)^\th q^{\th(\th-1)} B_q(\mu_i+(n-i)\th,\th) \,.
$$
Therefore the result equals
$$
C(\mu,n) H(\mu)
$$
times the following power of $q$
$$
\multline
\th |\mu| + \th^2 n(n-1)/2 -\th(\th-1)(n-1)\\
-\th(\th-1)(n-1)(n-2)/2
-\th(n-1)(n-2)/2-\th(n-1)\,,
\endmultline
$$
which equals $\th |\mu|$. \qed
\enddemo

Write
$$
\mu\subseteq\l
$$
if $\l_i\ge\mu_i$ for all $i$, and
$$
\mu\subset\l
$$
if $\mu\subseteq\l$ and $\mu\ne\l$\,. In the sequel
we shall deal with some sums like
$$
\sum_\nu c_\nu P_\nu
$$
where the only thing that matters is which $P_\nu$
can enter the sum, no matter with which specific constant factors
$c_\nu$. Constant factor like $c_\nu$ will always denote
such unspecific factors. 

The following proposition is a strengthening of a result by
S.~Sahi and F.~Knop.In particular, it gives the highest degree
term of $\Psm$, which was found in \cite{K,S}
by different methods.

\proclaim{Proposition 4.5}
$$
\multline
\Psm(x_1,x_2,\dots,x_n)=P_\mu(x_1,x_2 t^{-1},\dots,\xt{n})+\\
\sum_{\nu\subset\mu} c_\nu  P_\nu(x_1,x_2 t^{-1},\dots,\xt{n})\,, 
\endmultline \tag 4.11
$$
where $c_\nu$ are some constants
\footnote
{
Explicit formulas for these coefficients are given
by a particular case of the {\it binomial theorem\/}
for Macdonald polynomials, see the recent paper \cite{Ok3}
}
 which depend on $\mu$ and
$n$.
\endproclaim

\demo{Proof}
Induct on $n$ and $|\mu|$. Suppose $\mu_n>0$.  By \tht{4.4}
we have
$$
\Psm(x_1,\dots,x_n)=q^{|\mu^{-}|} 
\left(\sum_{k=0}^n (-t^{-n})^k \, e_{n-k} (\dots,\xt{i},\dots) \right)
P^*_{\mu^{-}}(x_1/q,\dots,x_n/q) \,. 
$$
By inductive assumption \tht{4.11} is true for $\mu^{-}$. It is well
known that the product
$$
e_r\, P_\nu
$$
is a linear combination of such Macdonald polynomials $P_\eta$
that
$$
\eta/\nu
$$
is a vertical $r$-strip. Therefore \tht{4.11} is true also for $\mu$.

Suppose $\mu_n=0$.  By inductive assumption
$$
\Psm(y_1,\dots,y_{n-1})=
P_\mu(y_1,\dots,y_{n-1}\,t^{2-n})+
\sum_{\nu\subset\mu} c_\nu  P_\nu(y_1,\dots,y_{n-1}\,t^{2-n})\,.
$$
Therefore
$$
\Psm(y_1,\dots,y_{n-1})=
t^{(2-n)|\mu|} P_\mu(\y_1,\dots,\y_{n-1})+
\sum_{\nu\subset\mu} c'_\nu  P_\nu(\y_1,\dots,\y_{n-1})\,,
$$
for some constants $c'_\nu$. Integrating as in \tht{3.2} and \tht{4.5}
we obtain
$$
\Psm(x_1,\dots,x_{n})=
t^{(1-n)|\mu|} P_\mu(\x_1,\dots,\x_n)+
\sum_{\nu\subset\mu} c''_\nu  P_\nu(\x_1,\dots,\x_{n})\,,
$$
which is equivalent to \tht{4.11}. \qed
\enddemo

 Using \tht{4.5} we
can reprove one more result from \cite{K},
which generalizes the corresponding result
for the shifted Schur functions \cite{OO}, Th.~3.1.
 
\proclaim{Proposition (Vanishing property) 4.6}
$$
\Psm(q^\l)=0,\quad\text{unless}\quad \mu\subseteq\l  \tag 4.12
$$
\endproclaim
\demo{Proof}
Induct on $n$ and $|\mu|$. If $\mu_n>0$ apply
proposition 4.2. 

If $\mu_n=0$ then apply the very same argument
which proved vanishing in the proof of the theorem.
\qed
\enddemo

\remark {\bf Remark 4.7} The formulas \tht{3.2} and \tht{4.5}
are certain $q$-integral operators
$$
\Lambda(n-1)\to\Lambda(n)
$$
and
$$
\Ls(n-1)\to\Ls(n)\,.
$$
Using projection \tht{3.12} (and its analog for
$\Ls$) in the
inverse direction we can characterize Macdonald
polynomials and shifted Macdonald polynomials
as eigenfunctions (with distinct eigenvalues) of 
some integral operators (in
the same way as for Schur functions \cite{OO},\S10).

In fact, there are countably many commuting integral
operators which correspond to iterated integration 
\tht{3.2} and \tht{4.5}
$$
\Lambda(n)\to\Lambda(N),\quad n<N\,,
$$
and projection back
$$
\Lambda(N)\to\Lambda(n)\,.
$$

\endremark

\head
\S5 Combinatorial formula for shifted Macdonald polynomials.
\endhead

In this section we shall establish branching rule
for shifted Macdonald polynomials, or (what is the
same) a combinatorial formula for $\Psm$ in terms
of semistandard tableaux on $\mu$.

Fist, we need some qualitative results about
the branching for $\Psm$. 

\proclaim{Proposition 5.1} Let $d$ be a variable, then
$$
\Psm(d x_1,\dots, d x_n) = \sum_{\nu\subseteq\mu}
c_\nu(d) P^*_\nu(x_1,\dots,x_n)\,, \tag 5.1
$$
where $c_\nu(d)$ are some polynomials in $d$ which
depend also on $\mu$ and $n$.
\endproclaim
\demo{Proof}
It is clear that
an expansion
$$
\Psm(d x_1,\dots, d x_n) = \sum_{\nu}
c_\nu(d) P^*_\nu(x_1,\dots,x_n)
$$
exists. We have to show that 
$$
c_\nu(d)=0, \quad\text{unless}\quad \nu\subseteq\mu \,.
$$
By \tht{4.11} the bases 
$$
\left\{\Psm(x_1,\dots,x_n)\right\} \tag 5.2
$$
and 
$$
\left\{P_\mu(\dots,\xt{i},\dots)\right\} \tag 5.3
$$
are mutually triangular with respect to the partial
ordering by inclusion of diagrams. Since dilatation
by $d$ in diagonal in the basis \tht{5.3} it is
triangular in the basis \tht{5.2}. \qed
\enddemo

Put
$$
\mu'=(\mu_1,\dots,\mu_{n-1})\,.
$$
We have
\proclaim{Proposition 5.2} 
$$
\Psm(x_1,\dots,x_{n-1},d) = \sum_{\nu\subseteq\mu'}
c_\nu(d) P^*_\nu(x_1,\dots,x_{n-1})\,, \tag 5.4
$$
where $c_\nu(d)$ are some polynomials in $d$ which
depend also on $\mu$ and $n$.
\endproclaim
\demo{Proof}
Again we have to show that certain summands cannot
occur.

If $d=1$ then by stability of $\Psm$
$$
\Psm(x_1,\dots,x_{n-1},1)=
\cases \Psm(x_1,\dots,x_{n-1}),&\mu_n=0\\
0,&\mu_n>0\,.
\endcases 
$$
For general $d$  introduce new variables $\xi_i$
$$
x_i= d\,\xi_i\,,\quad i=1,\dots,n-1\,,
$$ and use \tht{5.1}. \qed 
\enddemo

Put
$$
\pmu=(\mu_2,\dots,\mu_{n})\,.
$$
Recall that 
$$
\nu\prec\mu
$$
means that
$$
\mu_1\ge\nu_1\ge\mu_2\ge\dots\ge\nu_{n-1}\ge\mu_n\,,
$$
which is equivalent to
$$
\pmu\subseteq\nu\subseteq\mu'\,.
$$

\proclaim{Proposition 5.3}
$$
\Psm(d,x_2,\dots,x_{n}) = \sum_{\nu\prec\mu}
\f\, P^*_\nu(x_2,\dots,x_{n})\,,  \tag 5.5
$$
where $\f$ are some polynomials in $d$.
\endproclaim

Note that because of stability of shifted Macdonald
polynomials the polynomial $\f$ does {\it not}
depend on $n$.

\demo{Proof}
By the shifted symmetry 
$$
\Psm(d,x_2,\dots,x_n)=\Psm(x_2/t,\dots,x_n/t,t^{n-1}d)\,.
$$
Hence by \tht{5.4}
$$
\Psm(d,x_2,\dots,x_n) = \sum_{\nu\subseteq\mu'}
c_\nu(t^{n-1}d) P^*_\nu(x_2/t,\dots,x_{n}/t)\,.
$$
By \tht{5.1} 
$$
P^*_\nu(x_2/t,\dots,x_n/t) = \sum_{\eta\subseteq\nu}
c_\eta(1/t) P^*_\eta(x_1,\dots,x_n)\,.
$$
Therefore only summands with
$$
\nu\subseteq\mu'
$$
can occur in \tht{5.5}.

Now assume that there is a summand in \tht{5.5} with
$$
\pmu\not\subseteq\nu\,.
$$
We can choose $\nu$ minimal, that is in such a way
that $P^*_\eta$ with $\eta\subset\nu$ do not 
enter the sum \tht{5.5}. Then for all sufficiently
large $d$
$$
\Psm(d,\nu_1,\dots,\nu_{n-1})\ne 0
$$
for 
$$
P^*_\nu(\nu_1,\dots,\nu_{n-1})\ne 0
$$
and all other summands in \tht{5.5} vanish by
the vanishing property \tht{4.12}.
Clearly, this contradicts \tht{4.12} and makes our
assumption impossible. \qed
\enddemo

Now recall the corresponding result for 
Macdonald polynomials
$$
P_\mu(d,x_2,\dots,x_{n}) = \sum_{\nu\prec\mu}
\p d^{|\mu/\nu|}\, P_\nu(x_2,\dots,x_{n})\,,  \tag 5.6
$$
where $\p$ are certain nonzero rational functions
of $q$ and $t$ which can be found
in \cite{M}, section VI.7. Introduce the following notation
$$
\fp{d\,}{\mn}=\prod_{s\in\mu/\nu} ( d -q^{a'(s)} t^{-l'(s)} )\,,
$$
where the numbers $a'(s)$ and $l'(s)$ for a square $s\in\mu$
were defined in section 4. The number 
$q^{a'(s)} t^{-l'(s)}$ is the $q$-analog of the
{\it content} of the square $s$. If
$$
\nu=\emptyset,\quad \mu=(r)\,,
$$
then
$$
\fp{d\,}{\mn}=\fp{d\,}{r}
$$
Therefore the  number $\fp{d\,}{\mn}$
is a generalization of the factorial power $\fp{d\,}{r}$. 
The main result of this section is the following

\proclaim{Theorem III}
$$
\Psm(d,x_2,\dots,x_{n}) = \sum_{\nu\prec\mu}
\p\, t^{-|\nu|}\,\fp{d\,}{\mn} \, P^*_\nu(x_2,\dots,x_{n})\,,  \tag 5.7
$$
\endproclaim

It is clear that iterating this formula we obtain the
semistandard tableaux sum formula \tht{1.4} for polynomials $\Psm$.
We shall use induction on $n$. In fact we shall need
only the following corollary of this theorem

\proclaim{Corollary 5.4} For all $r=1,\dots,n$.
$$
\Psm(x_1,\dots,x_r,q^{\mu_{r+1}},\dots,q^{\mu_n}) =
c\, x_1^{\mu_1}\dots x_r^{\mu_r} +\dots\,,\tag 5.8
$$
where $c$ is a nonzero factor and dots stand for
lower monomials in lexicographic order.
\endproclaim
\demo{Proof}
Given a partition
$$
\eta=(\eta_1,\dots,\eta_{n-1})
$$
set
$$
\i\eta=(\eta_i,\dots,\eta_{n-1})\,.
$$
Then
$$
c=t^{-\sum_1^r |\ip\mu|} 
\left(\prod \psi_{\i\mu,\ip\mu}\right) H({}^{(r+1)}\mu) \,.\qed
$$
\enddemo

\demo{Proof of the theorem} 
Induction on $n$. The case $n=1$ is clear.

Suppose $n>1$. We shall find
$$
|\mu/\nu|
$$
distinct zeros of the polynomial $\f$. 
Fix some $i$ and show that  
$$
\f=0,\quad d=q^{\mu_i-1}t^{1-i},\dots,q^{\nu_i}t^{1-i} \,. \tag 5.9
$$
We shall prove \tht{5.9} by induction on $\i\nu$, in other
words we shall deduce \tht{5.9} from the assumption that
$$
\fe=0,\quad d=q^{\mu_i-1}t^{1-i},\dots,q^{\eta_i}t^{1-i} \tag 5.10
$$
for all $\eta$ such that
$$
\i\eta\subset\i\nu\,.
$$
Suppose $d$ is as \tht{5.9}. We have
$$
\Psm(d,x_2,\dots,x_i,q^{\nu_{i}},\dots,q^{\nu_{n-1}})=
\sum_{\eta\prec\mu}
\fe\,  P^*_\eta(x_2,\dots,x_i,q^{\nu_{i}},\dots,q^{\nu_{n-1}})\,. \tag 5.11
$$
By the vanishing property \tht{4.12} only summands satisfying
$$
\i\eta\subseteq\i\nu
$$
are nonzero. On the other hand, if
$$
\i\eta\subset\i\nu
$$
then in particular
$$
\eta_i\le\nu_i
$$
and by our assumption \tht{5.10} the corresponding
summand vanishes. Therefore only summands with
$$
\i\eta=\i\nu
$$
enter the sum. By \tht{5.8} each summand has the following 
form
$$
c \fe\, \big( x_2^{\eta_1}\dots x_i^{\eta_{i-1}} +\dots
\big)\,,\tag 5.12
$$
where $c$ is a nonzero factor and dots stand for
lower monomials in lexicographic order. 

On the other hand, by shifted symmetry
$$
\Psm(d,x_2,\dots,x_i,q^{\nu_{i}},\dots,q^{\nu_{n-1}})=
\Psm(x_2/t,\dots,x_i/t, t^{i-1} d,q^{\nu_{i}},\dots,q^{\nu_{n-1}})\,.\tag 5.14
$$
By the vanishing property \tht{4.12}
this should vanish if $d$ is as in \tht{5.9}
and
$$
x_i/t=q^{\l_i}\,,
$$
for all sufficiently large integers
$$
\l_2\ge\dots\ge\l_i \,.
$$
Therefore for such $d$ the polynomial \tht{5.11} should
be identically zero.
By virtue of \tht{5.12} it is impossible unless 
$$
\fe=0
$$
for all $\eta$ such that $\i\eta=\i\nu$. This proves \tht{5.9}.

Since 
$$
\deg\f\le|\mu/\nu|
$$
the polynomial $\f$ equals
$$
\fp{d\,}{\mn}
$$
up to a factor. This factor is clear from \tht{4.10}
and \tht{5.6}. \qed
\enddemo

\example{Example 5.5} The shifted analogs of the 
elementary symmetric functions are 
$$
P^*_{(1^k)}(x)=\sum_{i_1<\dots<i_k} t^{k-\sum i_s}
\, \prod_s (x_{i_s} - t^{s-k}) \,. \tag 5.14
$$
\endexample

\remark{\bf Remark 5.6} The shifted symmetry of $\Psm$
results in some complicated identities for
coefficients $\p$ which are not clear from
the explicit formula for these coefficients.
\endremark

\head
\S6 Duality
\endhead

The duality we shall discuss in this section
relates shifted Macdonald polynomials with
parameters $q$ and $t$ to shifted Macdonald
polynomials with parameters $1/t$ and $1/q$.
Denote by
$$
\Ls_t
$$
the algebra $\Ls$ constructed in section 4.
Denote by
$$
\pk(x;t)=\sum_i t^{k(1-i)} (x_i^k-1) \tag 6.1
$$
the power-sum generators of $\Ls_t$. Consider the
following isomorphism $\os$
$$
\align
\Ls_t &@>\quad\os\quad>>\Ls_{1/q},\\
\pk(t)\,&\mapsto\frac{1-\,q^k}{1-(1/t)^k}\,\, \pk(1/q) \,. \tag 6.2
\endalign
$$
Note a slight difference with the
Macdonald duality automorphism. The isomorphism \tht{6.2}
has the following clear combinatorial 
interpretation which generalizes the corresponding
result for shifted Schur function \cite{OO}, Th.~4.1.

\proclaim{Proposition 6.1} For all $f\in\Ls_t$
$$
\left[\os(f)\right](t^{-\l})= f(q^{\l'}) \,. \tag 6.3
$$
\endproclaim
\demo{Proof}
It suffices to check \tht{6.3} on generators \tht{6.1}.
For any partition $\l$ we have
$$
\frac{1}{q-1}\sum_i t^{1-i} (q^{\l_i}-1)=
\frac{1}{1/t-1}\sum_j q^{j-1} (t^{-\l'_j}-1)\,, \tag 6.4
$$
which is example VI.5.1(a) in \cite{M}. To see \tht{6.4}
one can use the same trick as in \cite{OO}, Th.~4.1,  and sum
$$
\sum_{(i,j)\in\l} q^{j-1} t^{1-i}
$$
first along rows and then along columns. It is clear
that \tht{6.4} means
$$
\frac{1-q}{1-1/t}\,\,p^*_{1} (t^{-\l};1/q)= p^*_{1} (q^{\l'};t) \,.
$$
Now replace $q$ by $q^k$ and $t$ by $t^k$ in \tht{6.4}
to obtain
$$
\frac{1-q^k}{1-(1/t)^k}\,\,p^*_{k} (t^{-\l};1/q)= p^*_{k} (q^{\l'};t) \,,
$$
as desired. \qed
\enddemo

\proclaim{Corollary 6.2} 
$$
\os(\Psm(q,t))= c \, \Psmp(1/t,1/q)\,, \tag 6.7
$$
where $c$ is a constant which depends on $\mu$.
\endproclaim
\demo{Proof} Follows from \tht{6.3} and definition of $\Psm$. \qed
\enddemo

Recall the following notation of Macdonald
$$
b_\l(q,t)=\prod_{s\in\l} \frac{1-q^{a(s)}t^{l(s)+1}}
{1-q^{a(s)+1}t^{l(s)}} \,.
$$
In particular, 
$$
b_{\l'}(t,q)=\prod_{s\in\l} \frac{1-q^{a(s)+1}t^{l(s)}}
{1-q^{a(s)}t^{l(s)+1}} \,.
$$

\proclaim{Theorem IV} 
$$
\os(\Psm(q,t))= (-t)^{|\mu|}b_{\mu'}(t,q)\,\Psmp(1/t,1/q)\,. \tag 6.8
$$
\endproclaim
\demo{Proof} 
It is clear that in order to calculate the constant in \tht{6.7}
we have to evaluate the equality \tht{6.7} at $t^{-\mu'}$. By
\tht{6.3} and definition of $\Psm$ we have
$$
\align
c&=\frac{H(\mu;q,t)}{H(\mu';1/t,1/q)}\\
&=\frac{t^{-2 n(\mu)} q^{n(\mu')} 
\prod_{s\in\mu} (q^{a(s)+1} t^{l(s)} -1) }
{q^{2 n(\mu')} t^{-n(\mu)} 
\prod_{s\in\mu} (q^{-a(s)} t^{-l(s)-1} -1) }\\
&=(-t)^{|\l|}b_{\l'}(t,q)\,,
\endalign
$$
because
$$
\sum_{s\in\mu} l(s) = n(\mu), \qquad  
\sum_{s\in\mu} a(s) = n(\mu')\,. \qed
$$
\enddemo

\example{Example 6.3}
Since the automorphism \tht{6.2} preserves the
degree we can look at the corresponding isomorphism
of the graded algebras
$$
\Lambda \to \Lambda\,,
$$
which differs from Macdonald automorphism
$$
p_k \to (-1)^{k-1} \frac{1-q^k}{1-t^k}\,p_k
$$
by factor
$$
(-t)^{\deg} \,,
$$
where $\deg$ is the degree of a polynomial.
Since the highest term of $\Psm$ is the
Macdonald polynomial $P_\mu$ and
$$
P_\mu(q,t)=P_\mu(1/q,1/t)\,,
$$
the theorem \tht{6.8} gives a new proof of the
duality theorem of Macdonald together with
computation of the constant factor.
\endexample

\head
\S 7 Degeneration
\endhead

In this section  we consider {\it shifted Jack polynomials}.
These polynomials were considered by F.~Knop and S.~Sahi in \cite{KS}.
These polynomials are indexed by partitions $\mu$, depend
on a parameter $\th$ and equal, by definition, to
$$
\Psm(x\,;\th)=\lim_{q\to 1} (q-1)^{-|\mu|}\,\Psm(q^x;q,q^\th)\,. \tag 7.1
$$
In this section let us change notations as follows.
The dependence on $q$ will be denoted explicitly
in all prelimit expressions like in the RHS of \tht{7.1}.
The expressions without $q$ will denote the corresponding
degenerations. For example, we write
$$
\fp{z\,;q}{r}=(z-1)(z-q)\dots(z-q^{r-1})\,. \tag 7.2
$$
This is the shifted Macdonald polynomial in one
variable
$$
P^*_r(z\,;q,t)=\fp{z\,;q}{r}\,. 
$$
The corresponding shifted Jack polynomial equals
$$
P^*_r(z\,;\th)=\fp{z}{r}\,,
$$
where
$$ 
\fp{z}{r}=z(z-1)\dots(z-r+1)\,. \tag 7.3
$$
In the same way write
$$
\fp{z}{\mn}=\prod_{s\in\mu/\nu} (z - a'(s)+\th\,l'(s))\,,
$$
this is a generalization of \tht{7.3}. 

The easiest way to see that \tht{7.1} is a polynomial of degree
$|\mu|$ is to look at the combinatorial formula for $\Psm(q,t)$.
The branching rule \tht{5.7} has the following limit
$$
\Psm(x_1,x_2,\dots,x_{n}\,;\th) = \sum_{\nu\prec\mu}
\p(\th)\,\fp{x_1}{\mn} \, P^*_\nu(x_2,\dots,x_{n}\,;\th)\,,  \tag 7.4
$$
where $\p(\th)$ are the  branching coefficients for ordinary
Jack polynomials, which are rational functions in $\th$
and can be found in \cite{St} and \cite{M}, section VI.10. Note that
$$
\th=1/\alpha\,,
$$
where $\alpha$ is the traditional parameter for Jack polynomials. 

Because of the symmetry and vanishing properties of 
the shifted Macdonald polynomials we have
$$
\align
\Psm(x\,;\th) &\text{  is symmetric in variables  } x_i-\th\,i\,, \tag 7.5\\
\Psm(\l\,;\th)&=0 \quad \text{unless} \quad \mu\subseteq\l\,, \tag 7.6\\
\Psm(\mu\,;\th)&=\prod_{s\in\mu} (a(s)+\th\,l(s)+1)\,. \tag 7.7
\endalign
$$
The product in \tht{7.7} is a $\th$-analog of the hook length
product; denote it by $H(\mu\,;\th)$.

Before discussing the limit of integral formula \tht{4.5}
for $\Psm(x\,;q,t)$ consider some elementary properties of
finite sums. We have
$$
\sum_{y=x_2}^{x_1} \fp{y}{r} = 
\frac{\fp{x_1+1}{r+1}-\fp{x_2}{r+1}}{r+1} \,, \tag 7.8
$$
provided $x_1-x_2$ is a nonnegative integer. Using \tht{7.8}
one can evaluate sums
$$
\sum_{y=x_2}^{x_1} f(y) \tag 7.9
$$
for arbitrary polynomials $f(y)$. If $x_1-x_2$ is not a nonnegative
integer then {\it define} \tht{7.9} using \tht{7.8}. In particular,
we have
$$
\sum_{y=x_2}^{x_1}=-\sum_{y=x_1+1}^{x_2-1} \,.
$$
Now suppose that  $f$ is a polynomial and
$$
\lim_{q\to1} \frac 1{(q-1)^{s}} f(q^y)= \fb(y) 
$$
for some $s$ and $\fb$. Then
$$
\lim_{q\to1} \frac 1{(q-1)^{s+1}} \int_{q^{x_2}}^{q^{x_1+1}}
f(y)\dq y= \sum_{y=x_2}^{x_1} \fb(y) \,, \tag 7.10
$$
which follows directly from definitions if $x_1-x_2$ is
an integer.

Now we calculate the limit of \tht{4.5}. Recall that
\cite{GR,1.10}
$$
\lim_{q\to 1} \G_q(z)=\G(z)\,,
$$
where $q$ tends to $1$ from below. Therefore
$$
\lim_{q\to 1} \frac{(q^z\ft}{(1-q)^{\th-1}} =
\frac{\G(z+\th-1)}{\G(z)} \,.
$$
In the sequel we assume that
$$
\th=1,2,\dots\,.
$$
In this case
$$
\frac{\G(z+\th-1)}{\G(z)}=z(z+1)\dots(z+\th-1)\,.
$$
Recall that we denote by $\Vs(x\,;q)$ the denominator
in \tht{4.5}
$$
\Vs(x\,;q)=\prod_{i<j} ( q^{\th(n-i)} x_i - q^{\th(n-j)} x_j)
\prod_{i\ne j} (q^{\th(j-i)+1} x_i/x_j\ft \,.
$$
We have
$$
\lim_{q\to 1} \frac{\Vs(q^x;q)}{(1-q)^{n(n-1)(\th+1/2)}}
= (-1)^{n(n-1)/2} \Vs(x)\,, \tag 7.11
$$
where
$$
\Vs(x)=
\prod_{i<j}(x_i-x_j+\th(j-i)) \prod_{i\ne j}
\fp{x_i-x_j+\th(j-i)+\th-1}{\th-1} \,. 
\tag 7.12
$$
Denote by
$$
\beta^*(y,x\,;q,q^\th)=
\prod_{i<j} (q^{\th(n-1-i)}y_i - q^{\th(n-1-j)}y_j) 
\prod_{i,j} (q^{\th(j-i)+1-\th}y_i/x_j\ft
$$
the density of the beta measure in \tht{4.5}. In the same way 
we have
$$
\lim_{q\to 1} 
\frac{\beta^*(q^y,q^x;q,q^\th)}
{(1-q)^{n(n-1)\th+(n-1)(n-2)/2}} =
(-1)^{(n-1)(n-2)/2}\beta^*(y,x\,;\th)\,, \tag 7.14
$$
where, by definition,
$$
\beta^*(y,x\,;\th)=
\prod_{i<j} (y_i - y_j +\th(j-i))
\prod_{i,j} 
\fp{y_i-x_j+\th(j-i)-1}{\th-1}\,.  \tag 7.15
$$
Using \tht{7.10}, \tht{7.11}, and \tht{7.14} we obtain the following
degeneration of \tht{4.5}
$$
\frac 1{\Vs(x)} \sum_{y\prec x} \, \Psm(y\,;\th) \,
\beta(y, x\,;\th)
=C(\mu,n) \, \Psm(x\,;\th) \,, \tag 7.16
$$
where
$$
\sum_{y\prec x} = \sum_{y_1=x_2}^{x_1} \dots 
\sum_{y_{n-1}=x_n}^{x_{n-1}}\,,
$$
and
$$
C(\mu,n)=\prod B(\mu_i+(n-i)\th,\th)\,.
$$
For $\th=1$ we obtain the coherence property of the
shifted symmetric functions \cite{OO}, \S10.

\enddocument

\end